\documentstyle[amsfonts,amssymb,12pt]{article}

\newcommand{\be}{\begin{equation}}
\newcommand{\ee}{\end{equation}}
\newcommand{\mb}{\Bbb  }
\newcommand{\beq}{\begin{eqnarray}}
\newcommand{\eeq}{\end{eqnarray}}
\newcommand{\beqst}{\begin{eqnarray*}}
\newcommand{\eeqst}{\end{eqnarray*}}

\newcommand{\ga}{\gamma}
\newcommand{\Ga}{\Gamma}

\newcommand{\de}{\delta}

\newcommand{\al}{\alpha}
\newcommand{\la}{\lambda}

\newcommand{\e}{\varepsilon}
\def\R{{\mathbb R}}

\def\C{{\mathbb C}}
\def\qed{{\hfill $\Box$}}
\newtheorem{theorem}{Theorem}[section]

\newtheorem{lemma}{Lemma}[section]

\title{ {\bf On deconvolution methods}}

\author{Alexander G. Ramm\\
LMA/CNRS, 31 Chemin J. Aiguier\\
Marseille 13402, cedex 20, France\\
Anahit Galstian\\
Department of Mathematics,Kansas State Univerisity,\\ 
Manhattan, KS 66506, USA \\ 
E-mail: ramm@math.ksu.edu\\
http://www.math.ksu.edu/\,$\widetilde{\ }$\,ramm}

\begin{document}
\thispagestyle{empty}

\date{}
\maketitle

\thispagestyle{empty}



\begin{abstract}
\noindent Several methods for solving efficiently the
one-dimensional  deconvolution problem are proposed.
The problem is to solve the Volterra equation
${\mathbf k} u:=\int_0^t k(t-s)u(s)ds=g(t),\quad 0\leq t\leq T$.
The data, $g(t)$, are noisy. Of special practical interest is the case
when the data
are noisy and known at a discrete set of times.
A general approach to the deconvolution problem is proposed:
represent ${\mathbf k}=A(I+S)$, where a method for a stable inversion of
$A$ is known, $S$ is a compact operator, and $I+S$ is injective.
This method is illustrated by examples: smooth kernels
$k(t)$, and weakly singular kernels, corresponding to Abel-type
of integral equations, are considered.
A recursive estimation scheme for solving deconvolution
problem with noisy discrete data is justified mathematically,
its convergence is proved, and error
estimates are obtained for the proposed deconvolution method.
\end{abstract}

\section{Introduction}

\setcounter{section}{1}
\setcounter{equation}{0}
\renewcommand{\theequation}{\thesection.\arabic{equation}} 

\bigskip

In many applications one models the relation between input signal $u(t)$ and
output signal $g(t)$ by the equation
\be\label{1.1}
{\mathbf k}u: = \int_0^t k(t-s)u(s) ds := k\ast u = g(t), \quad 0\le t\le 
T,
\ee
where $k(t)$, given for all $t\geq 0$, characterizes the linear system, 
$k\ast u$ is the 
convolution,
$u(t)=k(t)=0$ for $t<0$ and the cases $T<\infty$ and $T=\infty$ are both
of interest. In practice $g(t)$ is measured with some error, so $g_\de(t)$
is known, \,$\|g_\de(t)-g(t)\|\le \de$.\, The norm  we use 
is $L^2(0, \infty;e^{-2\sigma t})$ \,norm,
or $L^2(0, T)$, or
$L^{\infty}(0, T)$, and the 
case
$T<\infty$ can be reduced to the case $T=\infty$, as we show below.

If \,the \,operator \,${\mathbf k}$ \, in \,(\ref{1.1})\, is \,considered 
\, as \,an \,operator
on \,$L^\infty(0,T)$,\, and\\ ${\int_0^T |k(t)|dt < \infty}$,\, then 
${\mathbf k}$
is not boundedly invertible, so problem (\ref{1.1}) 
is ill-posed. One can see this   
from the formula \,${\int_0^t k(t-s)e^{ins}ds \to 0}$\,
as \, $n \to \infty $.

If $T<\infty$, one sets $u(s)=0$ for $t>T$ 
and defines $g(t)$ for $t>T$ as the
left-hand side of equation (\ref{1.1}). 
If this is done, then (\ref{1.1}) can be considered
as an equation on $(0,\infty)$ and its solution equals to $u(t)$,  the solution
of (\ref{1.1}) on $[0,T]$, when $t\in [0,T]$. We assume that 
$k(t)\not\equiv 0$
and \,${\int_0^\infty |k(t)|dt < \infty}$. Then  (\ref{1.1}) has at most
one solution in $L^1[0,T]$ if $k(t)\not= 0$ 
almost everywhere in $[0,T]$ (\cite{T}, p.327).
The assumption ${\int_0^\infty |k(t)|dt < \infty}$ can be replaced
without loss of generality
by a weaker assumption ${\int_0^\infty \exp(-\sigma t)|k(t)|dt < \infty}$,
with an arbitrary large fixed $\sigma >0$. This weaker assumption can be 
reduced to the original one by changing variables.

 A deconvolution method is a method to construct a stable approximation
$u_\de(t)$ of the solution $u(t)$ to equation (\ref{1.1}), given $g_\de$:
\be
\label{1.2}
\|u_\de(t)-u(t)\|:=\eta(\de)\to 0 \quad {\rm as}\quad \de\to 0. 
\ee
An operator $R(\de)$ which constructs such $u_\de$ from $g_\de$, 
$u_\de = R(\de)g_\de$, is called a regularizer (or a regularizing family, 
since $\de\to 0$) if (\ref{1.2}) holds.

There is a large literature on ill-posed problems. General methods
 for constructing
regularizers have been developed. They include variational regularization,
iterative regularization,  method of quasisolutions, etc [9]. 
In Section 2 the specific form of equation (\ref{1.1})
is used for constructing regularizers for equation (\ref{1.1}).
The emphasis is on the causality property of the regularizer.
The idea is similar to the one in [6] and [3]. In Section 
3 a simple general method to construct regularizers for equation (1.1)
is proposed. This method is practically efficient. It is  
illustrated by two examples in which the results from [2]-[5],  and [7]  
are used.
In Section 4 we investigate a recursive algorithm proposed in [1] for 
solving equation (\ref{1.1}) with noisy discrete data.
Again, the emphasis is on the causality property of the estimate: we use
only the data collected up to the time $t$ in order to estimate 
the signal $u$ at this time. If one would use a variational regularization
for constructing a stable estimation of $u(t)$, one would have to
use all the data collected on the full time inteval $[0,T]$, and not only 
on the "current" time interval $[0,t]$.  
Our analysis is much shorter than in [1] and yields more detailed results.
Moreover, In Section 6 we discuss briefly a generalization of these 
results to the case of
operator-valued kernels, which includes, in particular, matrix-valued 
kernels, that is, systems of Volterra equations. In 
Section 5 proofs are given.

\section {A deconvolution method}

\setcounter{equation}{0}
\renewcommand{\theequation}{\thesection.\arabic{equation}} 

\bigskip

Let $K(\la):= {\displaystyle \int_0^\infty e^{-\la t}k(t)dt}$. 
By  capital letters
the Laplace transform is denoted. If $u$ solves (\ref{1.1}) then 
$U(\la)=K^{-1}(\la)G(\la)$. One has:
\be
\label{2.1}
u_\de(t):= \frac{1}{2\pi i}\int_{\sigma -i\infty}^{\sigma + i\infty} e^{\la t}
K^{-1}(\la)G_\de(\la) 
\frac{1}{\left(\frac{\la}{N}+1\right)^m}d\la, \quad \la=\sigma + i\mu\in 
C_{\sigma},
\ee
where $C_{\sigma}$ is the contour $\Re \lambda=\sigma$, 
$\lambda=\sigma+i\mu$, $m>0$ is a 
sufficiently large positive 
integer and $N>0$ is a large parameter. We do not show the dependence of 
$u_\de(t)$ on $m$ and $N$ to simplify the notations.
We want to prove that one can choose $N=N(\de)\to \infty$ 
as $\de\to 0$, such that
(\ref{1.2}) holds. Note that $G_\de(\la)=G(\la)+W(\la)$, where 
${ |W(\la)|= \left| \int_0^\infty e^{-\la t}w(t)dt\right|\le 
\frac{\de}{\sigma}}$,
the noise $w(t)$ satisfies the inequality $|w(t)|\le \de, \sigma = \Re \la >0$,
and  $G(\la):=  \int_0^\infty e^{-\la t}g(t)dt$. 
When one uses $L^{\infty}(0, T)$ norm one assumes that
\beq
\label{2.2}
  K^{-1}(\la) G(\la)\in L^1(C_\sigma),
\eeq
This is an a priori assumption on $u(t)$. 

We assume throughout this paper that:
\be
\label{2.3}
|K(\la)|\ge c|\la|^{-a}\,\,, \quad \la\in C_\sigma, \,\,\, a\in {\mb R}.
\ee
Here and below $c>0$ denote various constants independent of $\de$ and $N$. The 
constant $a$ may be negative, but in many applications $a\geq 0$.

\medskip

\begin{theorem}
If $m>a+0.5$ and $(\ref{2.3})$ holds, then there exists 
$N(\de)\to \infty$\, as \, $\de \to 0$, such that $(\ref{1.2})$ holds
with $L^2(0,\infty; e^{-2\sigma t})$ norm.

If (2.2) and (2.3) hold, and $m>a+1$, then (1.2) holds with 
$L^\infty (0,T)$ norm.

If $m>a+1$, (2.3) holds, and
\be
\label{2.4}
|K^{-1}(\la)G(\la)| \le \frac{c}{1+|\la|^{1+d}}, \quad d={\rm const}>0\,,
\la \in C_\sigma,
\ee
then $N(\de)=O \left(\de^{-\frac{1}{q}}\right)$, $q=a+d+1$,  
 $\eta(\de) = O(\de^{\frac d q})$ if $0<d<1$,
and the norm in (1.2) is $L^\infty(0,T)$ norm.

If $d>1$, then $N(\de)=O(\de^{-\frac 1{2+a}})$ and $\eta(\de) = 
O(\de^{\frac 1 {a+2}}).$

If $d=1$, then $N(\de)=O(\exp[ \frac {|\log \de|}{a+2}(1+o(1)])$, and
$\eta(\de) = O(\frac {|\log \de|}{\exp[\frac  {|\log 
\de|}{a+2}(1+o(1)]}).$
\end{theorem}

We sketch proofs in the last Section.

\section{A general approach to deconvolution}

\setcounter{equation}{0}
\renewcommand{\theequation}{\thesection.\arabic{equation}} 

Suppose the operator ${\mathbf k}$ in (\ref{1.1}) can be decomposed into a sum
${\mathbf k}:= A+B$, where $A^{-1}B:= S$ is compact in the Banach space $X$, in
which ${\mathbf k}$ acts, and $I+S$ is boundedly invertible, or which is the 
same by the Fredholm alternative,
$N(I+S)= \{0\}$, where $N(A)$ is the null space of $A$. In this case $I+S$
is an isomorphism of $X$ onto $X$, ${\cal R}(A) = {\cal R}({\mathbf k})$, 
and
\be
\label{3.1}
{\mathbf k}u=A(I+S)u = g\,.
\ee
If a regularizer for $A$ is known, then (\ref{3.1}) can be solved stably 
by the
scheme
\be
\label{3.2}
u_\de = (I+S)^{-1}R(\de)g_\de \,\,,
\ee
and (\ref{1.2}) holds.

Since $I+S$ is an isomorphism, the error $\|v-v_\de\|$ of the approximation of
the stable solution of the equation $Av=g$ by the formula $v_\de= R(\de)g_\de$
is of the same order as $\|u_\de-u\|$.

\medskip

\noindent{\bf Example 3.1.}\, Let $k(t)\in C^1(0,T)$ and assume 
$k(0)\not=0$. 
Then,
without loss of generality, one may assume $k(0)=1$. Write (\ref{1.1})  as
\[
{\mathbf k}u= \int_0^t u(s)ds + \int_0^t \left[ k(t-s)-1\right] u(s) ds := Au + Bu = g.
\]
Here stable inversion of $A$ is equivalent to  stable numerical differentiation of
noisy data. This problem has been solved in \cite{R1} (see also \cite{R2} 
- \cite{R4}, [7]), and the results of these works yield the 
following theorem:

\begin{theorem} 
Assume that  $\|u\|_{W^{2,\infty}} \le 
m_2<\infty$. Then the operator
$$R(\de)g_\de:= \frac{g_\de(t+h(\de))-g(t)}{h(\de)}
$$ 
is a regularizer for the operator 
$A$ if $h(\de) = 2\left(\frac{\de}{m_2}\right)^{1/2}$.
One has: $\|R(\de)g_\de - u\|_{L^\infty_{[0,T-h]}} \le 2\sqrt{m_2\de}\,$. 
\end{theorem}
In [5] weaker a priori assumption on $u$ is used: $\|u\|_a \le m_a$,
$0<a\leq 1$, where the Hoelder-space norm is defined as 
$\|u\|_a:= \displaystyle{\sup_{x\not=y, x,y\in [0,T]}}
\frac{|u(t)-u(s)|}{|t-s|^a} + 
\displaystyle{\sup_{0\le  t \le T}} |u(t)|\,$.

In this example our method yields the equation 
$$(I+S)u_\de = R(\de)g_\de,$$
where $S$ is a Volterra operator: $Su_\de : = \int_0^t k'(t-s)u_\de(s)ds$.
 Therefore $u_\de$ can be easily found by iterations.

\medskip

\noindent {\bf Example 3.2.}\, Let $k(t)=\frac{t^{-\ga}}{\Ga (1-\ga)}+ 
m(t)$,\, $0<\ga<1$,\,$m(t)\in C^1$, 
$$Au=\frac{t^{-\ga}}{\Ga(1-\ga)}\ast u,\quad Bu=m\ast u.
$$
 One has 
$A^{-1}g= \frac{1}{\Ga(\ga)}\int_0^t \frac{g'(s)ds}{(t-s)^{1-\ga}}$.\,
Define 
$$R_1(\de)g_\de:=  \frac{1}{\Ga(\ga)}\int_0^t 
\frac{(R(\de)g_\de)(s)}{(t-s)^{1-\ga}}ds,$$
where $R(\de)$ is defined in Theorem 3.1. 

\begin{theorem}
 The operator-function $R_1(\de)$ is a regularizer for 
equation 
$Au=g$, the operator $S:= A^{-1}B$ is compact in $L^2(0,T)$ and 
$N(I+S)=\{0\}$.
\end{theorem}
Compactness of $S$ is clear from its definition. The operator
$R_1(\de)$ is a regularizer for $A$ because $R(\de)$ is a regularizer for the 
operator of differentiation. Finally $N(I+S)=\{0\}$ because
$A^{-1}B$ is a Volterra operator. Therefore $u_\de$ can be
easily computed by our general method and (\ref{1.2}) holds.

\section{Recursive estimation given discrete noisy data}

\setcounter{equation}{0}
\renewcommand{\theequation}{\thesection.\arabic{equation}} 

Assume  that 
$$g_\de(nh):= \xi_n = g(nh)+ w_n
$$ 
are noisy measurements of the data 
$g(t)$ at the time moments $nh$, $h>0$ is small number, $|w_n|\le \de$ is noise.
One wishes to estimate stably $u(t)$, the solution to $(1.1),$ given the 
data $\xi_1,....\xi_n$.
The following estimation method was essentially proposed in [1].
Set $v_0= \frac{\xi_1 - \xi_0}{h}$. Define $v_j$
by recursive formulas:

\be
\label{4.1}
\al v_n + \sum_{j=0}^{n-1}\int_{jh}^{(j+1)h}k(nh-s)ds  v_j=\xi_n, \,\quad 
 v_0= \frac{\xi_1-\xi_0}{h}\,.
\ee
In this Section we assume that $k(t)\in L^1(\R_+)\cap L^2_{loc}(\R_+)$
and $u(t)\in  L^2_{loc}(\R_+)$. Then $g(t)\in C_{loc}(\R_+)$, as
the following lemma claims.

{\bf Lemma 4.1}  {\it If $u,k\in L^2_{loc}(\R_+)$ then $k \ast u \in 
C_{loc}(\R_+)$.}

{\bf Proof.} One has
\beqst
|(k\ast u)(t+h)-(k\ast u)(t)|& \leq &  \int_0^{t+h}|k(t+h-s)-k(t-s)||u|ds\\
&+ &\int_t^{t+h}|k(t-s)||u|ds  := I_1+I_2.
\eeqst
Now,
$$I_1^2\leq \int_0^{t+h}|k(t+h-s)-k(t-s)|^2ds \int_0^{t+h}|u|^2ds \to 0 
 \hbox { as } h\to 0,
$$ 
and
$$I_2^2\leq \int_t^{t+h}|k(t-s)|^2ds  \int_t^{t+h}|u(s)|^2ds\to 0   
\hbox { as } h\to 0.
$$
Lemma 4.1 is proved. $\Box$

This lemma makes it reasonable to assume that $g(t)\in C_{loc}(\R_+)$.
In [1] the case $g\in L^2_{loc}(\R_+)$ is discussed, when $g$ is not
defined pointwise. It is proposed in [1] to use a mollification of $g$ 
around the points $nh$ instead of using $g(nh)$. However, this 
mollification 
requires a knowledge of $g$ in a neighborhoods of all points $nh$, and
this is an information different from the one assumed
at the beginning, namely  $\xi_1,....\xi_n$.
By this reason and because of Lemma 4.1, we assume that $g(t)\in 
C_{loc}(\R_+)$, the space of functions continuous on any compact subinterval
of $\R_+$.

Our assumptions in this section are:

{\bf A)} $k(t)\in L^1(\R_+)\cap L^2_{loc}(\R_+)$
and $u(t)\in  L^2_{loc}(\R_+)$,

{\bf B)} (2.3) holds,

{\bf C)} the union of the spectra of
 $\{K(\la)\}_{\forall \la\in C_\sigma}$ does not contain the set
$\{z: z\in \C,
\pi -\varphi < \arg z < \pi+\varphi, |z|<r\}\,,$
 where $\varphi > 0$ and $r>0$ are  arbitrary small fixed numbers.

If the assumptions {\bf A)}, {\bf B)}, and {\bf C)} hold, then our result 
is Theorem 4.1 
below. If, in addition, $k$ and $g$ are Hoelder-continuous,
then our result is Theorem 4.2 below, which gives the rate of convergence
in (1.2).
 
We prove that if $\al=\al(\de)$ and $h=h(\de)$ are chosen 
suitably, then the function $v_\de(t)$, defined by the formula:
$$
v_\de(t)=v_j \quad  for \quad jh\le t\le(j+1)h,
$$ approximates stably $u(t)$, so 
that 
$\|v_\de(t)-u(t)\|\to 0$
as $\de \to 0$, where $||\cdot||$ is $L^\infty(0,T)$-norm.
The rate of convergence is estimated in Theorem 4.2 under additional a 
priori assumptions. 

Let
\be
\label{4.2}
\al u_\al(t) + \int_0^t k(t-s)u_\al(s) ds = g(t), \quad 0\le t \le T.
\ee
The function $v_\de(t),$ defined above, solves the equation
\be
\label{4.3}
\al v_\de(t) + \int_0^t k(t-s)v_\de(s) ds = f_\de(t), 
\ee
where $f_\de(jh)=\xi_j$ for $t=jh,$ 
and for other values of $t$ the function $f_\de(t)$ is 
defined
as the left-hand side of (\ref{4.3}):
\[
f_\de(t):= g_\de(nh)-g(t) + \int_0^t \left[ k(t-s)-k(nh-s)\right]v_\de(s) ds 
+ g(t):= \varphi _\de(t) + g(t),
 \]
where $(n-1)h\leq t \leq nh, \quad n=1,2,3,......$.

One has 
$$|\varphi_\de(t)|\le \de+\ga_g(h) + c(T)\|v_\de\|\ga_k(h),$$ 
where
$$\ga_g(h):= |g(nh)-g(t)| \to 0,\quad  \ga_k(h):= 
\int_0^t |k(t-s)-k(nh-s)| ds\to 0 \quad \hbox { as } h\to 0,$$ 
and  \, $(n-1)h \le t\le nh,$ $n=1,2,3.....$. 

Denote $v_\de - u_\al:=w$. From (4.2) and (4.3) one gets:
$w=(\al+{\mathbf k})^{-1}\varphi_\de$, where $\al+{\mathbf k}:=
\al I+{\mathbf k}$, and $I$ is the identity operator.
One has $||(\al+{\mathbf k})^{-1}||\leq c \al^{-1}$, and
$||\varphi_\de||\leq \de + \ga_g(h)+c(T)\ga_k(h)\|v_\de\|.$ 
Therefore:
\[
\|v_\de - u_\al\|\le c(T) \frac{\de+\ga_g(h)+c(T)\ga_k(h)\|v_\de\|}{\al}.
\]
Let us choose 
$h=h(\de)$ such that 
$\frac{\ga_k(h)+\ga_g(h)}{\al}\to 0$ as $\al\to 0$, and $\al=\al(\de)$ 
so that $\frac{\de}{\al(\de)}\to 0$ as $\de \to 0$.

Let $u_\de:= u_{\al(\de)}$. Then (\ref{1.2}) holds by Lemma 5.1 below:
$$||u_\de-u||=||u_{\al(\de)}-u||=||(\al+{\mathbf k})^{-1}{\mathbf 
k}u-u||:=\varepsilon (\al) \to 0 \hbox { as } \al:=\al(\de)\to 0.$$ 

>From (\ref{4.3})
one gets 
\be
\label{4.4}
\hspace*{1.5cm} v_\de= (\al+{\mathbf k})^{-1}g + 
(\al+{\mathbf k})^{-1}[g_\de(nh)-g(t)] + (\al+{\mathbf k})^{-1}
{\mathbf Q} v_\de,
\ee
where 
$${\mathbf Q}v_\de:= \int_0^t [k(t-s)-k(nh-s)]v_\de ds.$$
 Since $\|(\al + {\mathbf k})^{-1}\|\le \frac{c}{\al}$
and $\|g_\de(nh)-g(t) \| \le c(T)(\de+\ga_g(h)),$
 one has: 
$$\|(\al + {\mathbf k})^{-1}[g_\de(nh)-g(t) ]\|
\le c(T)\frac{\de+\ga_g(h)}{\al}\to 0\quad  \de\to 0,$$
provided that $\al=\al(\de)$ and $h=h(\de)$ are chosen so that 
$\frac{\de}{\al(\de)}\to 0$ and 
$\frac{\ga_g(h(\de))}{\al(\de)}\to 0$  as $\de \to 0.$
Let $\al=\al(\de)$. Then
$\|(\al + {\mathbf k})^{-1}{\mathbf Q}v_\de \|\le 
\frac{c}{\al}\ga_k(h)\|v_\de\|,$ and if
$\frac{\ga_k(h(\de))}{\al(\de)}\to 0$ as $\de \to 0$, 
then $\|(\al + {\mathbf k})^{-1}{\mathbf Q}\|\leq 
\frac{c \ga_k(h)}{\al}\to 0$
as $\de \to 0$. Thus 
$$v_\de = [I-(\al + {\mathbf k})^{-1}{\mathbf Q}]^{-1}
\left[(\al + {\mathbf k})^{-1}g + O\left(\frac{\de+\ga_g(h)}{\al}\right)
\right],$$
where $O\left(\frac{\de+\ga_g(h)}{\al}\right)$ denotes an element
whose norm is  $O\left(\frac{\de+\ga_g(h)}{\al}\right)$. 
 Therefore, with
 $u_\de:=(\al(\de) + {\mathbf k})^{-1}g,$ one gets: 
$$\|v_\de-u_\de\|\le c(T) \frac{\de+\ga_g(h)+\ga_k(h)}{\al}.$$
Consequently:
\be
\label{4.5} 
\|v_\de - u\|\le \|v_\de-u_\de\| + \|u_\de -u\|\le c(T) \frac{\de+\ga_g(h)+\ga_k(h)}{\al}+ \varepsilon(\al)\to 0
\quad {\rm as }\quad \de\to 0.
\ee
We have proved:
\begin{theorem}
Assume {\bf A), B)}, and {\bf C)}. Then there exist $h=h(\de)\to 0$ and 
$\al=\al(\de)\to 0$
such that $(\ref{4.5})$ holds.
\end{theorem}

If $g(t)$ and $k(t)$ are Hoelder-continuous, 
then $\ga_g(h)+\ga_k(h)=O(h^b), \, 0<b\leq 1$.
Put $\de=h^b$. Then (4.5) can be written as:
\be
\label{4.6}
\|v_\de - u\|\le c[ \frac{\de}{\al}+\varepsilon (\al)].
\ee
Let us estimate $\varepsilon (\al)=||u_{\de}-u||$. One has:
\beqst
\varepsilon(\al)& = & \int_{-\infty}^\infty \left|\al(\al+K(\la))^{-1}\right| 
\left|K^{-1}(\la)G(\la)\right|d\mu\\
& \le & \int_{-\infty}^\infty \al|(\al+K(\la))^{-1}| \frac{d\mu}
{[1+(\sigma^2+\mu^2)^{1/2}]^{1+d}}\\
&= & 
 \int_{-M}^M  + \int_{-\infty}^M + \int_{M}^\infty :=
I_{1}+I_{2}+ I_{3}.
\eeqst
Using (\ref{2.3}), one gets: 
$$I_{1}\le c\al \int_{-M}^M  \frac{\mu^{a}d\mu}{(1+|\mu|)^{d+1}} 
\le c\al M^{a-d}.
$$
The estimates of \,$I_{2}$\, and  \,$I_{3}$\, are similar. Let us
estimate, for example, $I_3$: 
 $$I_{3} \le c \int_{M}^\infty  \frac{d\mu}{(1+\mu)^{d+1}}\leq c M^{-d}, 
\quad d>0.
$$
 Thus 
$$\varepsilon(\al)\le c(\al M^{a-d}+ 
M^{-d}).$$
If $d\geq a$, then
$\varepsilon(\al)\le c \al$, $||v_{\de}-u||\le c[ \frac{\de}{\al}+
\varepsilon (\al)].$  Minimizing with respect to $\al$, one gets
$||v_{\de}-u||\le c\de^{1/2}$ if $\al=\de^{1/2}.$ 
If $d<a$, then $\varepsilon (\al)\le c(\al M^{a-d} +M^{-d})$.
Minimizing with respect to $M$, one gets
$\varepsilon (\al)\le c\al^{d/a}$, $||v_{\de}-u||\le c(\de \al^{-1}+ 
\al^{d/a})$. Minimizing with respect to $\al$, one gets
$||v_{\de}-u||\le c\de^{\frac {d}{d+a}}.$  

Let us summarize the result:

\begin{theorem}
Assume {\bf A), B)}, and {\bf C)}. If
$g$ and $k$ are Hoelder-continuous, so that $\gamma_g (h)
+\gamma_k (h)=O(h^b)$,
$0<b\leq 1$, then $||v_{\de}-u||\le c\de^{\frac {d}{d+a}},$
provided that $h=\de^{\frac 1 b}$, $\al=O(\de^{\frac a {d+a}})$, and 
$d<a$.
If $d\geq a$, then $||v_{\de}-u||\le c\de^{\frac {1}{2}}.$
\end{theorem}

In \cite{R7} a singular perturbation
problem was solved for a class of one- and multidimensional integral
equations. The problem we study in Sec. 4 contains a singular perturbation 
problem as a basic component: we are interested in the behavior of the
operator $(\al + {\mathbf k})^{-1}$ as $\al \to +0$.

\section{Proofs.}

\setcounter{equation}{0}
\renewcommand{\theequation}{\thesection.\arabic{equation}} 

The norm below is $L^2(0,\infty; e^{-2\sigma t})$ norm, it is equivalent 
to $L^2(0,T)$ norm on $(0,T)$.
By the spectrum of a scalar function $K(\la)$ we mean the set of its
values, and if $K(\la)$ is an operator-valued function, then its spectrum
is defined as usual.

\begin{lemma}
Let (2.3) hold and assume that the union of the spectra of 
 $\{K(\la)\}_{\forall \la\in C_\sigma}$ does not contain 
the set 
$
\{z: z\in \C, 
\pi -\varphi < \arg z < \pi+\varphi, |z|<r\}\,,$
 where $\varphi > 0$ and $r>0$
are  arbitrary small fixed numbers. Then 
\[\varepsilon (\al):=
\|(\al + {\mathbf k})^{-1}{\mathbf k}u -u\|\to 0 \quad as \quad\al\to +0
\quad for \, all \quad u \in L^2(0,\infty; e^{-2\sigma t})\,.
\]
\end{lemma}

\medskip

\noindent {\bf Proof.} One has 
\beqst
& & \|(\al + {\mathbf k})^{-1}{\mathbf k}u -{\mathbf k}^{-1}{\mathbf k}u\|^2 = 
\|(\al + {\mathbf k})^{-1}\al{\mathbf k}^{-1}{\mathbf k}u\|^2\\
& & =   \|\al (\al + {\mathbf k})^{-1}u\|^2 =
\left\|\frac{e^{\sigma t}}{2\pi}\int_{-\infty}^\infty e^{i\mu t}\al(\al+K(\la))^{-1}
 K^{-1}(\la)G(\la)d\mu \right\|^2\\
& & = \frac{1}{2\pi}\int_{-\infty}^\infty 
\al^2\left|(\al+K(\la))^{-1}\right|^2 
\left|K^{-1}(\la)G(\la) \right|^2d\mu:= 
\e(\al)\to 0 \quad {\rm as} \quad \al\to 0\,,\\
& &  \la= \sigma+i\mu \in C_\sigma.
\eeqst
Here we have used: 1) Parseval's equality; 2) the asumption $u\in 
L^2(0,\infty; 
e^{-2\sigma t})$ which is equivalent to $K^{-1}(\la)G(\la) \in L^2(C_\sigma)$;
3) the estimate $\displaystyle {\sup_{\la\in 
C_\sigma}}\bigg|\al(\al+K(\la))^{-1}\bigg|$ $\le$  $c$
 which follows
from the assumption about the range of $K(\la)$ on 
$C_\sigma$;  4) estimate (2.3); and 
5) the dominated convergence theorem. Lemma 5.1 is proved  \qed

\begin{lemma}
Under the assumptions of Lemma 5.1 one has:
\[
\|(\al + {\mathbf k})^{-1}g_\de -u\| \le \frac{c\de}{\al} +\e(\al), 
\quad \e(\al)\to 0 \quad {\rm as} \quad \al\to 0, \quad c=const>0\,. 
\]
\end{lemma}

\medskip

\noindent {\bf Proof.} One has 
\[
\|(\al + {\mathbf k})^{-1}g_\de -u\| \le 
\|(\al + {\mathbf k})^{-1}(g_\de -g)\| +
\|(\al + {\mathbf k})^{-1}g -u\| 
\le \|(\al + {\mathbf k})^{-1})\| \de +\e(\al)\,. 
\]
By Lemma 5.1, $\e(\al)\to 0$ as $\al\to 0$, and 
$\|(\al + {\mathbf k})^{-1}\|\le c\al^{-1}$,
where $c=$const$>0$ depends on $\varphi$, 
as follows from the proof of Lemma 5.1 and from the estimate 
$$
 \sup_{\la\in C_\sigma}|(\al+K(\la))^{-1}|\le \frac {c}{\al}.
$$ 
Lemma 5.2 is proved. $\Box$

\medskip

\noindent {\bf Proof of Theorem 2.1.} If $u_\de$ is defined in 
(\ref{2.1}) and $\|\cdot\|$
is $L^2(0,\infty; e^{-2\sigma t})$ norm, then Parseval's equality yields
\beqst
& &\|u_\de-u\|^2 =  \frac{1}{2\pi}\int_{-\infty}^\infty \left|K^{-1}(\la)G_\de(\la)
\frac{1}{\left(\frac{\la}{N}+1\right)^m} - K^{-1}(\la)G(\la)\right|^2d\mu\\
& &\le 
\frac{1}{\pi}\int_{-\infty}^\infty \left|K^{-1}(\la)G(\la)\right|^2
\left|\frac{1}{\left(\frac{\la}{N}+1\right)^m} - 1\right|^2d\mu 
 + \frac{\de^2}{\pi\sigma^2} \int_{-\infty}^\infty \frac{|K^{-1}(\la)|^2}{ 
\left|\frac{\la}{N}+1\right|^{2m}}d\mu\\
& &:= I_1 + \de^2 I_2 \,,
\eeqst
where we have used the formulas $G_\de=G+W$, $|W|\leq \frac {\de}{\sigma}.$

If $N\to\infty$ then $I_1=I_1(N)\to 0$ by the dominated convergence 
theorem.
Let us estimate $I_2$ assuming (\ref{2.3}) and taking $m-a>0.5$:
\beqst
I_2 &\le & c N^{2m}\int_{-\infty}^\infty 
\frac{(\sigma^2+\mu^2)^a d\mu }{[(\sigma + N)^2+\mu^2]^m}
\le  c N^{2m}\int_{-\infty}^\infty \frac{d\mu }{[(\sigma + N)^2+\mu^2]^{m-a}}\\
&\le & c \frac{N^{2m}}{(N+\sigma)^{2m-2a-1}}
\int_{-\infty}^\infty \frac{d\nu }{(1+ \nu^2)^{m-a}}
\le c N^{2a+1}\,.
\eeqst
Thus, if $2a+1 > 0$, then, using the estimate
$(x+y)^{1/2}\leq x^{1/2}+y^{1/2},$ $x,y \geq 0,$  one gets:
\[
\|u_\de-u\|\le \de c^{1/2} N^{a+\frac{1}{2}} + I^{\frac{1}{2}}_1(N):= 
\eta(\de,N)\,.
\]
Minimizing $\eta(\de,N)$ with respect to $N$
for a fixed $\de$, denoting the minimizer by $N(\de)$, $N(\de)\to \infty$,
as $\de\to 0$, and the minimum by $\eta(\de):= \eta(\de,N(\de)),$ 
one gets
$\eta(\de)\to  0$ as $\de\to 0$. Thus, (\ref{1.2}) is proved 
with $L^2(0,\infty; e^{-2\sigma t})$ norm. $\Box$

If $L^\infty(0,T)$ norm is used for $\|u(t)-u_\de(t)\|$, then one gets 
(\ref{1.2}) 
if assumptions (\ref{2.2}) and  (\ref{2.3}) are used.
Namely, 
\beqst
\sup_{0\leq t \leq T}|u_{\de}(t)-u(t)|& \le & 
\frac{\exp(\sigma T)}{2\pi}\int_{-\infty}^{\infty}
|K^{-1}G_\de (1+\frac {\la}{N})^{-m}-K^{-1}G|d\mu\\
& \le & c(T)(J_1+J_2):=\eta,
\eeqst
where
$$J_1:=J_1(N)=\int_{-\infty}^{\infty}|K^{-1}G||(1+\frac {\la}{N})^{-m}-1|d\mu,$$
and
$$J_2:=J_2(N)=\frac {\de}{\sigma}\int_{-\infty}^{\infty}|K^{-1}||1+\frac 
{\la}{N}|^{-m}d\mu.$$
If $K^{-1}G\in L^1(C_\sigma)$, then $\lim_{N\to \infty}J_1=0$ by the 
dominated convergence theorem.
If $m>a+1$ and (2.3) holds, then 
$$J_2\leq c\de N^{a+1}.$$ 

{\it Thus (1.2) holds with
$L^\infty(0,T)$-norm provided that 
$K^{-1}G\in L^1(C_\sigma)$, $m>a+1$ and (2.3) holds.}

{\it If (2.4) holds, one can get a rate of decay.}
Namely
$$J_1=\int_{-M}^M+\int_{|\mu|>M}=j_1 +j_2,$$
where 
\beqst
j_1 & \leq & \int_{-M}^M (1+|\la|^{d+1})^{-1}|(1+\frac {\la}{N})^{-m}-1|d\mu\\
& \leq & c \int_{0}^M (1+|\la|^{d+1})^{-1}\frac {|\la|}{N}d\mu\leq 
cM^{1-d}N^{-1}, \hbox { if }  0<d<1,
\eeqst
$$
j_1\leq cN^{-1} \hbox { if }  d>1,
$$
$$
j_1\leq c \frac {\log M} {N}  \hbox { if }  d=1,
$$ 
and 
$$ j_2\leq c\int_{M}^\infty (1+|\la|^{d+1})^{-1}|(1+\frac 
{\la}{N})^{-m}-1|d\mu \leq  c \int_{M}^\infty (1+|\la|^{d+1})^{-1} 
\leq cM^{-d},  
$$
$c>0$ stands for various constants, and $\la=\sigma +i\mu$. Thus
$$J_1\le c(M^{1-d}N^{-1}+M^{-d})  \hbox { if }  0<d<1, 
$$
$$J_1\le c(N^{-1}+M^{-d})  \hbox { if }   d>1,   
$$
$$J_1\le c(\frac {\log M} {N}+M^{-d})  \hbox { if }  d=1.   
$$
If $0<d<1$ then choose $M=N$ and get $J_1\le cN^{-d}$.
Therefore, if $m>a+1$ and (2.3) and (2.4) hold, then
$$J_1+J_2\le c(\de N^{a+1}+N^{-d})    \hbox { if }  0<d<1.$$
Minimizing with respect to N for a fixed $\de>0$, one
gets the minimizer $N=N(\de)=O(\de^{-\frac 1{1+a+d}})$
and the estimate $\eta \le O(\de^{\frac {d}{1+a+d}})$ if $ 0<d<1$.
If $d>1$ then choose $M^d=N$ and get $J_1 \leq cN^{-1}$,
$J_1+J_2\le c(\de N^{a+1}+N^{-1})$, $N(\de)=O(\de^{-\frac 1{2+a}})$,
and the estimate $\eta \le O(\de^{\frac {1}{2+a}})$.
If $d=1$ then choose $M=N$ and get $J_1 \leq c \frac {\log N} N$,
$J_1+J_2\le c(\de N^{a+1}+\frac {\log N} N)$,
$N(\de)=O(\exp[ \frac {|\log \de|}{a+2}(1+o(1)])$,
and the estimate $\eta \le O(\frac {|\log \de |}{\exp[\frac 
{|\log \de |}{2+a}(1+o(1))]})$.
Theorem 2.1 is proved.  $\Box$

\medskip

\noindent{\bf Remark 5.3.} If $k(t)>0$ and $k(t)\to 0$ monotonically, then
$\Re K(\la)\ge 0$ for $\Re \la >0$. If $\Re K(\la)\ge 0$ then 
$\bigg|\al(\al+K(\la))^{-1}\bigg|\le 1$. Condition $\Re K(\la)\ge 0$ implies
that the assumption of Lemma 5.1 holds with $\varphi=\frac{\pi}{2}$
and $r=\infty$.

\section{Generalizations.}
Most of our results and proofs remain valid
for operator-valued functions $k(t)$, in particular
for matrix-valued kernels, that is, for systems of
Volterra equations. Let $k(t)$ be an
operator in a Banach space, and $K(\la)$ be 
its Laplace transform. If one replaces
the absolute values by the norms in (2.3),
(2.4) and elsewhere in the proofs, then one gets 
Theorems 2.1, 3.1, 3.2, and 4.1
and lemmas 5.1 and 5.2 with operator-valued $k(t)$.

\medskip

{\bf Acknowledgement.} AGR thanks Prof. L. Pandolfi and Dr. F. Fagnani 
for discussions. Their paper [1] was useful for the results in Sec. 4 
of this paper.

\end{document}